\documentclass[leqno]{amsart}
\usepackage{amsmath}
\usepackage{amssymb}
\usepackage{amsthm}
\usepackage{enumerate}
\usepackage[mathscr]{eucal}
\theoremstyle{plain}
\usepackage{tikz}
\newtheorem{theorem}{Theorem}[section]

\newtheorem{prop}[theorem]{Proposition}
\theoremstyle{definition}
\newtheorem{definition}[theorem]{Definition}
\newtheorem{remark}[theorem]{Remark}

\newtheorem{cor}[theorem]{Corollary}
\theoremstyle{remark}

\begin{document}
	
	\title [On norm derivatives and the ball-covering property of Banach spaces] {On norm derivatives and the ball-covering property of Banach spaces} 
	
	\author{Debmalya Sain}

	\address{(Sain)~Department of Mathematics, Indian Institute of Information Technology Raichur, Karnataka 584135, India}
	\email{saindebmalya@gmail.com}

	%

	\thanks{The author feels elated to acknowledge Professor Miguel Martin and Dr. Saikat Roy for their kind encouragement while preparing this manuscript.}
	
	
	

	\date{}
	
	\begin{abstract}
		We study a local version of the ball-covering problem in Banach spaces, and obtain a complete solution to it in terms of the norm derivatives. We illustrate the advantage of the local approach by obtaining substantial refinements of several previously known results on this topic.
	\end{abstract}
	
	\maketitle 
	
	\section{Introduction}
	The unit ball of a Banach space largely determines its algebraic, geometric, and topological properties. Indeed, the problem of understanding the characteristic features of a Banach space essentially reduces to studying the structure of the unit ball of the space. Motivated by this simple observation, several mathematicians have investigated the properties of the unit ball of a given Banach space, from different perspectives. Cheng initiated the study of the ball-covering property of a Banach space in \cite{C}, and demonstrated that many important properties of a Banach space, such as smoothness, have deep connections with the ball-covering properties of the space. In light of the seminal work done in \cite{C}, many important applications of the ball-covering property have been obtained in later studies, and the topic remains active till date. We refer the readers to \cite{CCS, CKWZ, CZZ, GLM, LLLZ}, and the references therein, for more information in this regard. The purpose of this article is to consider a generalized version of the ball-covering property considered in \cite{C}, and to study the corresponding problem \emph{from a purely local point of view.}\\
	
	The letter $ X $ denotes a Banach space over the real field, and let $ B_X, S_X $ denote the unit ball and the unit sphere of $ X $, respectively. A non-zero element $ x \in X $ is said to be a smooth point in $ X $ if the set $ J(x) := \{ f \in S_{X^{*}} : f(x) = \| x \| \} $ is a singleton. Geometrically, a unit vector $ x \in S_X $ is smooth in $ X $ if there exists a unique supporting hyperplane to $ B_{X} $ at $ x. $ The open ball centered at $ x $ and having radius $ r > 0 $ is denoted by $ B(x, r) := \{ y \in X : \| x-y \| < r \}. $ Throughout this article, we will be considering open balls without further mentioning it explicitly. In \cite{C}, Cheng considers the following\\
	
	\noindent \textbf{Problem 1:} Given a Banach space $ X, $ how many (in the sense of cardinality) balls not containing the zero vector can together cover the unit sphere $ S_X $ of $ X? $ In particular, if $ X $ is finite-dimensional, is there a smallest number of balls having such property?\\
	
	\noindent In this article, We will be focusing on a local version of Problem $ 1 $ by considering the following\\ 
	
	\noindent \textbf{Problem 2:} Given non-zero elements $ x, y $ in a Banach space $ X, $ does there exist a ball centered at $ \lambda x, $ for some $ \lambda \in \mathbb{R}, $ which contains $ y $ but does not contain the zero vector?\\
	
	\noindent We obtain a complete solution to Problem $ 2 $ and further show that the local nature of it allows us to answer Problem $ 1 $ from a much broader perspective. In particular, this local approach to the ball-covering problem will enable us to replace $ S_X $ in Problem $ 1 $ by an arbitrary bounded set $ A \subset X $ which is at a positive distance from the zero vector and satisfies the condition that given any non-zero $ x \in X, $ there exists $ r_x > 0 $ such that $ r_x x \in A. $\\
	
	The answer to Problem $ 2 $ is obtained in terms of the so called norm derivatives which play an important role in understanding the geometry of Banach spaces.  
	
	\begin{definition}
		Let $ \mathbb{X} $ be a real normed linear space and let $ x, y \in \mathbb{X}. $ The norm derivatives at $ x $ in the direction of $ y $ are defined as 
		\[\rho'_{+} (x, y) = \lim_{t \rightarrow 0+ } \|x\| \frac{\|x+ty\| - \|x\|}{t},\] 
		
		\[\rho'_{-} (x, y) = \lim_{t \rightarrow 0- } \| x \| \frac{\|x+ty\| - \|x\|}{t}.\]
		
	\end{definition} 
	 
	 Norm derivatives have been studied extensively by many mathematicians due to its applicability in several geometric problems in the setting of Banach spaces. The readers are referred to \cite{AST, CW, KS, S, Wa} for some interesting applications of the norm derivatives and related properties of Banach spaces, including smoothness. We next state, mainly for the convenience of the readers, some of the important facts regarding the norm derivatives, which find applications in our present work. Since all these results are well-known, we do not provide the proofs, and refer the readers to \cite{AST} for a detailed treatment of the same.\\
	 
	 \begin{itemize}
	 	\item The convexity of the norm function ensures that the mappings $ \rho'_{\pm} $ are well-defined.\\
	 	
	 	\item Given any $ x, y \in X $ and any $ \alpha \in \mathbb{R}, $ the following statements hold true:
	 	\noindent \[ (i)~ \rho'_\pm(\alpha x, y)=\rho'_\pm(x,\alpha y) = \left\{
	 	\begin{array}{ll}
	 		\alpha \rho'_\pm(x,y),  & \mbox{if } \alpha\geq 0 \\
	 		\alpha \rho'_\mp(x,y),  & \mbox{if } \alpha< 0.
	 	\end{array}
	 	\right.\]

	 	$ (ii)~ \rho'_{-} (x, y) \leq \rho'_{+} (x, y).~ \textit{Moreover, a non-zero $ x \in X $ is smooth in $ X $ if and only if}~ \rho'_{+} (x, y) = \rho'_{-} (x, y) $ for all $ y \in X. $ \\
	 	
	 	$ (iii)~ \rho'_{\pm} (x, \alpha x + y) = \alpha \| x \|^2 + \rho'_{\pm} (x, y). $ \\

	 	$ (iv)~ \rho'_{+} (x, y) = \|x\| \sup \{ x^*(y) : x^* \in \mathbb{J}(x) \}. $ \\
	 	
	 	$ (v)~ \rho'_{-} (x, y) = \|x\| \inf \{ x^*(y) : x^* \in \mathbb{J}(x) \}. $ \\
	 	
	 	$ (vi)~ $ Given any $ x, y \in X, $ $ \| x + \lambda y \| \geq \| x \| $ for all scalars $ \lambda $ if and only if $ \rho'_{-} (x, y) \leq 0 \leq \rho'_{+} (x, y). $ In this context, we recall from \cite{B, J} that given any $ x, y \in X, $ we say that $ x $ is Birkhoff-James orthogonal to $ y, $ written as $ x \perp_B y, $ if $ \| x + \lambda y \| \geq \| x \| $ for all scalars $ \lambda. $ Moreover, the James characterization of Birkhoff-James orthogonality \cite{J} states that $ x \perp_B y $ if and only if there exists $ f \in J(x) $ such that $ f(y) = 0. $ We use the notation $ x^{\perp} $ to denote the set of all vectors $ y $ such that $ x \perp_B y. $\\
	
	\end{itemize}
	
As we will see in the next section, a direct and complete answer to Problem $ 2 $ is given by the sign of the norm derivative $ \rho'_{\pm} (x, y). $ Moreover, this local approach will allow us to obtain refinements of several well-known results related to the ball-covering properties of a Banach space.

	\section{Main Results}

	We begin with two central results of this article, which together completely characterize the local version of the ball-covering property in Banach spaces in terms of the norm derivatives.
	
	\begin{theorem}\label{theorem:norm derivatives and ball-covering1}
		Let $ X $ be a Banach space and let $ x, y \in X $ be non-zero. Then the following two conditions are equivalent:\\
		$ (i) $ There exists a ball centered at $ \lambda x, $ for some $ \lambda > 0, $ which contains $ y $ but does not contain the zero vector.\\
		$ (ii) $ $ \rho'_{-} (x, y) > 0. $
    \end{theorem}

	\begin{proof}
	By virtue of the positive homogeneity property of $ \rho'_{-}, $ we may and do assume without any loss of generality that $ \| x \| = 1. $ We establish the theorem by proving that $ (i) \iff (ii) $ in the following two steps:\\
	
	\noindent $ (i) \implies (ii): $ Let $ B(\lambda x, r) $ be the ball centered at $ \lambda x $ that contains $ y $ but does not contain the zero vector. We note that $ 0 < r \leq \lambda. $ Since $ \| \lambda x - y \| < r \leq \| \lambda x \|, $ it follows that $ y \notin (\lambda x)^{-}. $ As $ \lambda > 0, $ by virtue of Proposition $ 2.1 $ of \cite{Sa}, this is equivalent to $ y \in x^{+} \setminus x^{-}. $ Therefore, $ f(y) > 0 $ for every $ f \in J(x). $ Since $ \rho'_{-} (x, y) = \| x \| \inf\{{f(y) : f \in J(x)}\}, $ we conclude that $ \rho'_{-} (x, y) \geq 0. $ We next claim that $ \rho'_{-} (x, y) > 0. $ Suppose on the contrary that $ \rho'_{-} (x, y) = 0. $ Then $ x \perp_B y, $ which implies that $ \| \lambda x - y \| \geq \| \lambda x \| = \lambda, $ a contradiction to our hypothesis $ y \in B(\lambda x, r), $ as $ r \leq \lambda. $ This establishes our claim.\\ 
	\noindent $ (ii) \implies (i): $ Clearly, it suffices to show that there exist $ \lambda > 0 $ and $ 0 < r_{\lambda} \leq \lambda $ such that $ y \in B(\lambda x, r_{\lambda}). $ Suppose on the contrary that $ y \notin \bigcup\limits_{\lambda > 0} B(\lambda x, r_{\lambda}), $ whenever $ 0 < r_{\lambda} \leq \lambda. $ We can write $ y = \beta x + z, $ for some $ \beta \in \mathbb{R} $ and some $ z \in X. $ Then for each $ \lambda > \max\{1, |\beta|\}, $ considering $ r_{\lambda} = \lambda - \frac{1}{\lambda} \in (0, \lambda), $ we obtain that
	\begin{align*}
		\lambda - \frac{1}{\lambda} \leq \| \lambda x - y \| &\implies -\frac{1}{\lambda} \leq \| \lambda x - \beta x -z \| - \lambda \\\
		&\implies -\frac{1}{\lambda} \leq \| (\lambda - \beta) x - z \| - \lambda \| x \| \\\
		&\implies -\frac{1}{\lambda} \leq (\lambda - \beta) [\| x - \frac{1}{\lambda - \beta} z \| - \| x \|] - \beta
		\end{align*}
	Taking $ t = \frac{1}{\lambda - \beta}, $ the last inequality can be rewritten in the following form:
	\[ -\frac{1}{\lambda} \leq \frac{\| x - tz \| - \| x \|}{t} - \beta. \]
	It is clear that $ t \longrightarrow 0+ $ as $ \lambda \longrightarrow +\infty. $ Therefore, letting $ \lambda \longrightarrow +\infty, $ we conclude that 
	\[ 0 \leq \rho'_{+} (x, -z) - \beta. \]
	By using the properties of norm derivatives, as mentioned in the introduction, and putting $ z = y - \beta x, $ we obtain the following chain of implications:
	\[ \rho'_{-} (x, -\beta x + y) \leq - \beta \implies - \beta \| x \|^2 + \rho'_{-} (x, y) \leq - \beta \implies \rho'_{-} (x, y) \leq 0, \]
	a contradiction to our hypothesis that $ \rho'_{-} (x, y) > 0. $ This completes the proof of the theorem.
	\end{proof}

We next characterize the counterpart of the local ball-covering property in Banach spaces, corresponding to the case $ \lambda < 0. $

\begin{theorem}\label{theorem:norm derivatives and ball-covering2}
	Let $ X $ be a Banach space and let $ x, y \in X $ be non-zero. Then the following two conditions are equivalent:\\
	$ (i) $ There exists a ball centered at $ \lambda x, $ for some $ \lambda < 0, $ which contains $ y $ but does not contain the zero vector.\\
	$ (ii) $ $ \rho'_{+} (x, y) < 0. $
\end{theorem}

    \begin{proof}
	We observe that by virtue of Theorem \ref{theorem:norm derivatives and ball-covering1}, the Condition $ (i) $ is equivalent to $ \rho'_{-} (-x, y) > 0. $ This in turn is equivalent to $ \rho'_{+} (x, y) < 0, $ due to the well-known properties of the norm derivatives. Hence the theorem.
	\end{proof}

Theorem \ref{theorem:norm derivatives and ball-covering1} and Theorem \ref{theorem:norm derivatives and ball-covering2} are of fundamental importance in studying the ball-covering properties of Banach spaces. We would like to emphasize that the true strength of these two results lie in their local nature, which allows us to obtain refinements of several well-known global results in this topic, besides some new observations. The following corollary, which follows directly from the above results, also gives a characterization of the deeply studied concept of Birkhoff-James orthogonality in Banach spaces in terms of ball-covering.

\begin{cor}\label{corollary:corollary1}
	Let $ X $ be a Banach space and let $ x, y \in X $ be non-zero. Then at most one of the following two conditions holds true:\\
	$ (i) $ There exists a ball centered at $ \lambda x, $ for some $ \lambda > 0, $ which contains $ y $ but does not contain the zero vector.\\
	$ (ii) $ There exists a ball centered at $ \lambda x, $ for some $ \lambda < 0, $ which contains $ y $ but does not contain the zero vector.\\
	Moreover, neither of the above two conditions holds true if and only if $ x \perp_B y. $  
\end{cor}

\begin{proof}
	Since $ \rho'_{-} (x, y) \leq \rho'_{+} (x, y), $ the first part follows trivially from Theorem \ref{theorem:norm derivatives and ball-covering1} and Theorem \ref{theorem:norm derivatives and ball-covering2}. Thereafter, the second part can be deduced directly from the well-known equivalence: 
	\[ x \perp_B y \iff \rho'_{-} (x, y) \leq 0 \leq \rho'_{+} (x, y). \]
\end{proof}

We next prove that with respect to the ball-covering property, balls centered on the same ray emanating from the zero vector can be replaced by a single ball having its center on the same ray.

\begin{prop}\label{proposition:proposition1}
	Let $ X $ be a Banach space and let $ A \subset X $ be such that $ A \subset \bigcup\limits_{i=1}^{n} B(\lambda_i x, r_i), $ where $ n \in \mathbb{N} $ and $ x \in S_{X} $ are fixed, $ \lambda_i > 0, $ and $ 0 < r_i \leq \lambda_i $ for each $ 1 \leq i \leq n. $ Then for each $ \lambda \geq \max\{ \lambda_i : 1 \leq i \leq n \}, $ there exists $ 0 < r_{\lambda} \leq \lambda $ such that $ A \subset B(\lambda x, r_{\lambda}). $
\end{prop}

\begin{proof}
	For a fixed $ \lambda \geq \max\{ \lambda_i : 1 \leq i \leq n \}, $ we set $ \mu_i = \lambda - \lambda_i \geq 0. $ Applying the triangle inequality for the norm function, it is easy to verify that $ A \subset \bigcup\limits_{i=1}^{n} B(\lambda_i x, r_i) \subset \bigcup\limits_{i=1}^{n} B(\lambda x, r_i + \mu_i). $ Moreover, for each $ 1 \leq i \leq n, $ it is trivially true that $ 0 < r_i + \mu_i \leq \lambda. $ Choosing $ r_{\lambda} = \max\{ r_i + \mu_i : 1 \leq i \leq n \}, $ we obtain that $ A \subset B(\lambda x, r_{\lambda}), $ thereby finishing the proof.
\end{proof}

In Theorem $ 2.2 $ of the pioneering article \cite{C}, Cheng has shown that for every $ n $-dimensional Banach space $ X, $ the unit sphere $ S_X $ has a symmetric ball-covering consisting of $ 2n $ balls. We obtain a refinement of this very interesting result for compact sets not containing the zero vector, in a Banach space which is not necessarily finite-dimensional. We recall that for a non-empty bounded set $ C $ of $ X, $ a point $ x_0 \in C $ is called an exposed point of $ C $ if there exists an $ x^{*} \in X^{*} $ such that 
\[ x^{*}(x_0) = \sup \{ x^{*}(x) : x \in C \}~~ \textit{and}~~ \{ x \in C : x^{*}(x) = x^{*}(x_0) \} = \{ x_0 \}. \]

In this case, we say that $ x^{*} $ is an exposing functional for $ x_0. $ When $ f_0 \in C \subset X^{*} $ and the exposing functional for $ f_0 $ is of the form $ \Psi(x), $ where $ x \in X $ and $ \Psi : X \longrightarrow X^{**} $ denotes the canonical embedding, we say that $ \Psi(x) $ is an weak*-exposing functional for $ f_0 \in C $ and that $ f_0 $ is an weak*-exposed point of $ C. $

\begin{theorem}\label{theorem:symmetric ball-covering}
	Let $ X $ be a Banach space and let $ A \subset X $ be a compact set in $ X $ not containing the zero vector. Let $ \{ f_i : 1 \leq i \leq m \} \subset S_{X^*} $ be a collection of weak*-exposed points of $ B_{X^{*}} $ such that $ \max\{ f_i(y) : 1 \leq i \leq m \} > 0 $ for each $ y \in A. $ Then $ A $ has a ball-covering consisting of $ m $ balls.
\end{theorem}

\begin{proof}
	Let $ \Psi(x_i) \in X^{**} $ be an weak*-exposing functional for $ f_i \in B_{X^{*}}, $ where $ x_i \in X $ and $ 1 \leq i \leq m. $ We claim that $ J(x_i) = \{ f_i \}, $ for each $ 1 \leq i \leq m. $ Indeed, it is clear that $ f_i \in J(x_i), $ since $ f_i(x_i) = \Psi(x_i)(f_i) = \| \Psi(x_i) \| = \| x_i \|. $ Moreover, if $ g \in J(x_i) $ then $ \Psi(x_i) (g) = g(x_i) = \| x_i \| = \Psi(x_i) (f_i), $ contradicting that $ \Psi(x_i) $ is the weak*-exposing functional for $ f_i. $ This establishes our claim and shows that each $ x_i $ is a smooth point in $ X. $ Given any $ y \in A, $ there exists $ f_y \in \{ f_i : 1 \leq i \leq m \} $ such that $ f_y(y) > 0. $ Let $ x_y \in X $ be such that $ \Psi(x_y) $ is the weak*-exposing functional for $ f_y. $ Since $ x_y $ is smooth in $ X, $ it follows that $ \rho'_{+} (x_y, y) = \rho'_{-} (x_y, y) = f_y(y) > 0. $ Applying Theorem \ref{theorem:norm derivatives and ball-covering1}, we conclude that there exists a ball $ B(\lambda_y x_y, r_y), $ where $ 0 < r_y \leq \lambda_y \| x_y \|, $ that contains $ y $ but does not contain the zero vector. Clearly, $ \{ B(\lambda_y x_y, r_y) : y \in A \} $ is an open cover of $ A. $ By using the compactness of $ A, $ we obtain a finite sub-cover $ \{ B(\lambda_{y_{k}} x_{y_{k}}, r_{y_{k}}) : y_k \in A, 1 \leq k \leq l \} $ of the open cover $ \{ B(\lambda_y x_y, r_y) : y \in A \}, $ where $ l \in \mathbb{N} $ is fixed. We note that $ x_y \in \{ x_i : 1 \leq i \leq m \} $ for each $ y \in A $ and moreover, $ \lambda_{y_{k}} > 0 $ for each $ 1 \leq k \leq l. $ Therefore, by virtue of Proposition \ref{proposition:proposition1}, for each $ 1 \leq i \leq m, $ there exist $ \lambda_i > 0 $ and $ 0 < r_i \leq \lambda_i \| x_i \| $ such that $ A \subset \bigcup\limits_{i=1}^{m} B(\lambda_i x_i, r_i). $ This completes the proof of the theorem.
\end{proof}

In case $ X $ is a finite-dimensional Banach space, Theorem \ref{theorem:symmetric ball-covering} admits the following substantial strengthening:

\begin{cor}\label{corollary:finite-dimensional1}
	Let $ X $ be a finite-dimensional Banach space and let $ A \subset X $ be a bounded set in $ X $ such that $ d(0, A) > 0. $ Let $ \{ f_i : 1 \leq i \leq m \} \subset S_{X^*} $ be a collection of weak* exposed functionals in $ X^* $ such that $ \max\{ f_i(y) : 1 \leq i \leq m \} > 0 $ for each $ y \in A. $ Then $ A $ has a ball-covering consisting of $ m $ balls.
\end{cor}

\begin{proof}
	The proof follows trivially from Theorem \ref{theorem:symmetric ball-covering}, in light of the fact that under these assumptions, $ \bar{A} $ is a compact set not containing the zero vector.
\end{proof}

\begin{cor}\label{corollary:corollary2}(Theorem $ 2.2 $ of \cite{C})
	Suppose that $ X $ is an $ n $-dimensional Banach space. Then:\\
	$ (i) $ $ S_X $ has a symmetric ball-covering consisting of $ 2n $ balls.\\
	$ (ii) $ Every symmetric ball-covering of $ S_X $ contains at least $ 2n $ balls.
\end{cor}

\begin{proof}
	$ (i): $ Since the set of exposed points of $ B_{X^*} $ is dense in the set of extreme points of $ B_{X^*}, $ it is easy to see from the Krein-Milman Theorem that there exists $ n $ linearly independent exposed points of $ B_{X^*}, $ say $ f_1, f_2, \ldots, f_n. $ In particular, it follows that $ \bigcap\limits_{i=1}^{n} ker f_i = \{ 0 \}. $ Therefore, given any $ x \in S_{X}, $ we get that $ \max\{ \pm f_i(x) : 1 \leq i \leq n \} > 0. $ Since $ S_X $ is a compact subset of $ X $ not containing the zero vector, the result follows from Theorem \ref{theorem:symmetric ball-covering}.\\
	$ (ii): $ If possible, suppose that $ \{ B(\pm x_i, r_i) : 1 \leq i \leq r \} $ is a symmetric ball-covering of $ S_X, $ where $ x_i \in X \setminus \{ 0 \}, $ $ 0 < r_i \leq \| x_i \|, $ and $ 1 \leq r < n. $ For each $ 1 \leq i \leq r, $ consider $ f_i \in J(x_i). $ Since $ 1 \leq r < n, $ there exists a unit norm vector $ z \in \bigcap\limits_{i=1}^{r} ker f_i. $ By the James characterization of Birkhoff-James orthogonality, we obtain that $ \pm x_i \perp_B z, $ for each $ 1 \leq i \leq r. $ Therefore, Corollary \ref{corollary:corollary1} asserts that $ z \notin \bigcup\limits_{i=1}^{r} \{ B(\pm x_i, r_i), $ a contradiction to our assumption that $ \{ B(\pm x_i, r_i) : 1 \leq i \leq r \} $ is a ball-covering of $ S_X. $ This proves that every symmetric ball-covering of $ S_X $ contains at least $ 2n $ balls.
\end{proof}

In Theorem $ 2.3 $ of \cite{C}, it is shown that whenever $ X $ is an $ n $-dimensional Banach space, $ S_X $ does not admit a ball-covering consisting of at most $ n $ balls. Moreover, it is also proved in the same theorem that the conclusions of Corollary \ref{corollary:corollary2} can be substantially improved, under the additional assumption that $ X $ is smooth. We are next going to present a refinement of this, which will allow us to deduce the original result as a direct consequence.

\begin{theorem}\label{theorem:general ball-covering}
	Let $ X $ be an $ n $-dimensional Banach space and let $ A \subset X $ be a bounded set in $ X $ not containing the zero vector. Also assume that given any non-zero $ x \in X, $ there exists $ r_x > 0 $ such that $ r_x x \in A. $ Then:\\
	$ (i) $ Every ball-covering of $ A $ contains at least $ n+1 $ balls.\\
	If, in addition, $ X $ is smooth and $ d(0, A) > 0, $ then\\
	$ (ii) $ $ A $ admits a ball-covering consisting of $ n+1 $ balls.
\end{theorem}

\begin{proof}
	$ (i): $ Suppose on the contrary that for some $ m \leq n, $ $ \{ B(\lambda_i x_i, r_i) : 1 \leq i \leq m \} $ is a ball-covering of $ A, $ where $ x_i \in X \setminus \{ 0 \}, $  and $ 0 < r_i \leq | \lambda_i | \| x_i \|. $ Without any loss of generality, we may and do assume that $ m = n, $ $ \lambda_i > 0, $ and $ \| x_i \| = 1, $ for each $ 1 \leq i \leq n. $ For $ 1 \leq i \leq n-1, $ let us choose $ f_i \in J(x_i) $ and keep it fixed throughout. Since $ X $ is $ n $-dimensional, we can find a unit norm vector $ z \in \bigcap\limits_{i=1}^{n-1} \ker f_i. $ Let $ z_0 = r_z z \in A, $ where $ r_z > 0. $ Since $ x_i \perp_B z $ for all $ 1 \leq i \leq n-1, $ using the homogeneity property of Birkhoff-James orthogonality, we obtain from Corollary \ref{corollary:corollary1} that $ z_0 \notin \bigcup\limits_{i=1}^{n-1} B(\lambda_i x_i, r_i). $ Therefore, $ z_0 \in B(\lambda_n x_n, r_n). $ Applying Theorem \ref{theorem:norm derivatives and ball-covering1}, we conclude that $ \rho'_{-} (x_n, z_0) > 0. $ Since $ \rho'_{+} (x_n, z_0) \geq \rho'_{-} (x_n, z_0), $ this implies that $ \rho'_{-} (x_n, -z_0) = -\rho'_{+} (x_n, z_0) < 0. $ Let $ w_0 = r_{-z_{0}} (-z_0) \in A. $ Clearly, $ \rho'_{-} (x_n, w_0) < 0. $ Moreover, it can be shown as before that $ w_0 \notin \bigcup\limits_{i=1}^{n-1} B(\lambda_i x_i, r_i). $ Therefore, $ w_0 \in B(\lambda_n x_n, r_n). $ By Theorem \ref{theorem:norm derivatives and ball-covering1}, this is equivalent to $ \rho'_{-} (x_n, w_0) > 0, $ a contradiction. This finishes the proof of $ (i). $\\
	$ (ii): $ Let $ \{ f_i  : 1 \leq i \leq n \} \subset S_{X^*} $ be any linearly independent set in $ X^* $ and let $ f_{n+1} = -\sum_{i=1}^{n} f_i. $ We claim that $ \max\{ f_i(y) : 1 \leq i \leq n+1 \} > 0 $ for each $ y \in A. $ If $ f_{i}(y) > 0 $ for some $ 1 \leq i \leq n, $ then we are done. So let us assume that $ f_{i}(y) \leq 0 $ for all $ 1 \leq i \leq n. $ Since $ 0 \notin A, $ and $ \{ f_i  : 1 \leq i \leq n \} $ is linearly independent, it follows in particular that there exists $ 1 \leq i_0 \leq n $ such that $ f_{i_{0}}(y) < 0. $ Then it is trivial to see that $ f_{n+1}(y) > 0, $ thereby justifying our claim. we now apply Corollary \ref{corollary:finite-dimensional1} to conclude that $ A $ admits a ball-covering consisting of $ n+1 $ balls. This completes the proof of $ (ii) $ and establishes the theorem.
\end{proof}

\begin{cor}\label{corollary:corollary3}(Theorem $ 2.2 $ of \cite{C})
	Suppose that $ X $ is an $ n $-dimensional Banach space. Then:\\
	$ (i) $ Every ball-covering of $ S_X $ contains at least $ n+1 $ balls.\\
	If, in addition, $ X $ is smooth, then\\
	$ (ii) $ $ S_X $ admits a ball-covering consisting of $ n+1 $ balls.
\end{cor}

\begin{proof}
	The proof follows trivially from Theorem \ref{theorem:general ball-covering}, by taking $ A = S_X. $
\end{proof}

In Theorem $ 2.6 $ of \cite{CCS}, the authors have characterized finite ball-coverings of $ S_X $ in terms of the subdifferential mapping. We recall from \cite{CCS} that the subdifferential mapping $ \partial \| . \| : X \longrightarrow 2^{B_{X^{*}}} $ of the norm is defined by $ \partial \| x \| := \{ x^* \in S_{X^{*}}: x^*(x) = \| x \| \}. $ Of course, in our terminology, $ \partial \| x \| = J(x) $ for every non-zero $ x \in X. $ We now obtain a refinement of Theorem $ 2.6 $ of \cite{CCS}, by replacing $ S_{X} $ with an arbitrary compact set not containing the zero vector.

\begin{theorem}\label{theorem:subdifferential}
	Let $ X $ be a Banach space and let $ A \subset X $ be a compact set in $ X $ not containing the zero vector. Suppose that $ I $ is an index set with $ m $ elements and $ \{ x_i : i \in I \} \subset S_{X}. $ Then $ \mathcal{B} \equiv \{ B(y_i, r_i) \}_{i \in I} $ forms a ball-covering of $ A $ for some $ y_i \in \mathbb{R}^+ x_i $ with $ \| y_i \| \geq r_i $ for all $ i \in I $ if and only if for every selection $ \phi $ of the subdifferential mapping $ \partial \| . \|, $ $ \{\phi(x_i)\}_{i \in I} $ positively separates points of $ A, $ that is, $ \sup_{i \in I} \phi(x_i) (x) > 0 $ for every $ x \in A. $ 
\end{theorem}

\begin{proof}
	Let us first prove the sufficient part of the theorem. Let $ x \in A $ be arbitrary but fixed after choice. We claim that there exists $ i_0 \in I $ such that $ \rho'_{-} (x_{i_{0}}, x) > 0. $ Suppose on the contrary that $ \rho'_{-} (x_i, x) \leq 0 $ for every $ i \in I. $ Then given any $ i \in I, $ there exists $ f_i \in J(x_i) $ such that $ f_i(x) \leq 0. $ Consider the selection $ \phi $ of the subdifferential mapping $ \partial \| . \| $ such that $ \phi (x_i) = f_i $ for all $ i \in I. $ Clearly, $ \sup_{i \in I} \phi(x_i) (x) \leq 0, $ a contradiction to our hypothesis. This proves our claim. Applying Theorem \ref{theorem:norm derivatives and ball-covering1}, we conclude that $ x \in B(\lambda_{i_{0}} x_{i_{0}}, r_{i_{0}}), $ where $ 0 < r_{i_{0}} \leq \lambda_{i_{0}}. $ Since $ x \in A $ was chosen arbitrarily, we can use the compactness of $ A $ in the same way as in the proof of Theorem \ref{theorem:symmetric ball-covering}, and then apply Proposition \ref{proposition:proposition1} to finish the proof.\\
	Let us next prove the necessary part of the theorem. Let $ \phi $ be any selection of the subdifferential mapping $ \partial \| . \| $ and let $ x \in A $ be arbitrary. Since $ \{ B(y_i, r_i) \}_{i \in I} $ forms a ball-covering of $ A, $ there exists $ i_0 \in I $ such that $ x \in B(\lambda_{i_{0}} x_{i_{0}}, r_{i_{0}}), $ where $ 0 < r_{i_{0}} \leq \lambda_{i_{0}}.  $ Once again applying Theorem \ref{theorem:norm derivatives and ball-covering1}, we get that $ \rho'_{-} (x_{i_{0}}, x) > 0. $ Using the expression for $ \rho'_{-} (x_{i_{0}}, x), $ it is clear that $ f(x) > 0 $ for all $ f \in J(x_{i_{0}}). $ In particular, this implies that $ \sup_{i \in I} \phi(x_i) (x) > 0, $ as desired. This completes the proof of the theorem.
\end{proof}

When $ X $ is a finite-dimensional Banach space, the necessary part of Theorem \ref{theorem:subdifferential} can be strengthened for bounded sets in the following way. 

\begin{theorem}\label{theorem:finite-dimensional subdifferential sufficient part}
	Let $ X $ be a finite-dimensional Banach space and let $ A \subset X $ be a bounded set in $ X $ such that $ d(0, A) > 0. $ Suppose that $ I $ is an index set with $ m $ elements and $ \{ x_i : i \in I \} \subset S_{X}. $ Also assume that $ \bar{A} \bigcap \bigcup\limits_{i \in I} x_{i}^{\perp} = \emptyset. $ If $ \mathcal{B} \equiv \{ B(y_i, r_i) \}_{i \in I} $ forms a ball-covering of $ A $ for some $ y_i \in \mathbb{R}^+ x_i $ with $ \| y_i \| \geq r_i $ for all $ i \in I, $ then for every selection $ \phi $ of the subdifferential mapping $ \partial \| . \|, $ $ \{\phi(x_i)\}_{i \in I} $ positively separates points of $ \bar{A}, $ that is, $ \sup_{i \in I} \phi(x_i) (x) > 0 $ for every $ x \in \bar{A}. $ 
\end{theorem}

\begin{proof}
	Let $ \phi $ be any selection of the subdifferential mapping $ \partial \| . \| $ and let $ x \in \bar{A} $ be arbitrary. There exists a sequence $ \{ a_n \} \subset A $ such that $ a_n \longrightarrow x $ as $ n \longrightarrow \infty. $ Applying Theorem \ref{theorem:norm derivatives and ball-covering1}, we obtain that for each $ n \in \mathbb{N}, $ there exists $ i_n \in I $ such that $ \rho'_{-} (x_{i_{n}}, a_n) > 0.  $ Since $ I $ is finite, we may and do assume without any loss of generality (by passing onto a subsequence, if necessary) that there exists $ i_0 \in I $ such that $ \rho'_{-} (x_{i_{0}}, a_n) > 0 $ for all $ n \in \mathbb{N}. $ Using the properties of $ \rho'_{-}, $ it is easy to see that $ \rho'_{-} (x_{i_{0}}, x) \geq 0. $ Also, $ \rho'_{-} (x_{i_{0}}, x) = 0 $ implies that $ x_{i_{0}} \perp_B x, $ a contradiction to our hypothesis. Therefore, we have the following chain of implications:
	\[ \rho'_{-} (x_{i_{0}}, x) > 0 \implies f(x) > 0 ~\forall f \in J(x_{i_{0}}) \implies \sup_{i \in I} \phi(x_i) (x) > 0. \]
	Since $ x \in \bar{A} $ is arbitrary, this finishes the proof.
\end{proof}

\begin{remark}
	Let $ A \subset X $ be a bounded set such that either of the following two conditions holds true:\\
	$ (i) $ $ d(0, A) = 0,~ (ii):~ \bar{A} \bigcap \limits_{i \in I} x_{i}^{\perp} \neq \emptyset.$\\
	Then it is not difficult to see that the conclusion of Theorem \ref{theorem:finite-dimensional subdifferential sufficient part} can no longer hold true. This illustrates that the conditions assumed in the previous theorem cannot be completely removed.
\end{remark}

On the other hand, it is possible to strengthen the sufficient part of Theorem \ref{theorem:subdifferential} for bounded sets, without any additional restrictions. In order to avoid the repetition of arguments, we only present a sketch of the proof.

\begin{theorem}\label{theorem:finite-dimensional subdifferential necessary part}
	Let $ X $ be a finite-dimensional Banach space and let $ A \subset X $ be a bounded set in $ X $ such that $ d(0, A) > 0. $ Suppose that $ I $ is an index set with $ m $ elements and $ \{ x_i : i \in I \} \subset S_{X}. $ If for every selection $ \phi $ of the subdifferential mapping $ \partial \| . \|, $ $ \{\phi(x_i)\}_{i \in I} $ positively separates points of $ \bar{A}, $ that is, $ \sup_{i \in I} \phi(x_i) (x) > 0 $ for every $ x \in \bar{A}, $ then $ \mathcal{B} \equiv \{ B(y_i, r_i) \}_{i \in I} $ forms a ball-covering of $ A $ for some $ y_i \in \mathbb{R}^+ x_i $ with $ \| y_i \| \geq r_i $ for all $ i \in I $  
\end{theorem}

\begin{proof}
	As in the proof of the sufficient part of Theorem \ref{theorem:subdifferential}, it can be deduced that given any $ u \in \bar{A}, $ there exists $ i_{u} \in I $ such that $ \rho'_{-}(x_{i_{u}}, u) > 0. $ Once again, we can use the compactness of $ \bar{A} $ in conjunction with Theorem \ref{theorem:norm derivatives and ball-covering1} and Proposition \ref{proposition:proposition1} to obtain the desired conclusion.
\end{proof}

In Proposition $ 2.7 $ of \cite{CCS}, the authors have established the importance of smooth points in the study of finite ball-coverings of the unit sphere of a separable Banach space. As the final result of this article, we obtain a strengthening of this useful observation by applying Theorem \ref{theorem:norm derivatives and ball-covering1}. Compatible with our main theme, we study a local version of this result by replacing separability with a much weaker condition of existence of a convergent sequence of smooth points. Of course, the corresponding global result, i.e., Proposition $ 2.7 $ of \cite{CCS} follows directly from the local result obtained by us.

\begin{theorem}\label{theorem:smoothness}
	Let $ x, y $ be non-zero elements of a Banach space $ X $ such that $ y \in B(\lambda x, r), $ where $ \lambda > 0 $ and $ 0 < r \leq \| x \|. $ Suppose that there exists a sequence $ \{ x_n \} $ of smooth points in $ X $ with $ x_n \longrightarrow x. $ Then there exists $ n_0 \in \mathbb{N} $ such that $ y \in B(\lambda_{n_{0}} x_{n_{0}}, r_{n_{0}}), $ where $ \lambda_{n_{0}} > 0 $ and $ 0 < r_{n_{0}} \leq \| x_{n_{0}} \|. $
\end{theorem}

\begin{proof}
	Without any loss of generality, we may and do assume that $ \| x \| = \| x_n \| = 1 $ for all $ n \in \mathbb{N}. $ Since $ y \in B(\lambda x, r), $ where $ \lambda > 0 $ and $ 0 < r \leq 1, $ it follows from Theorem \ref{theorem:norm derivatives and ball-covering1} that $ \rho'_{-} (x, y) > 0. $ Moreover, by virtue of the same result, it suffices to show that $ \rho'_{-} (x_{n_{0}}, y) > 0, $ for some $ n_0 \in \mathbb{N}. $ As each $ x_n $ is a smooth point in $ X, $ let us assume that $ J(x_n) = \{ f_n \} $ for all $ n \in \mathbb{N}. $ We now complete the proof by considering the two possible cases:\\
	Case I: $ \{ f_n : n \in \mathbb{N} \} $ is finite. Then there exists $ n_0 \in \mathbb{N} $ such that $ J(x_n) = \{ f_{n_{0}} \} $ for infinitely many $ n. $ Passing onto a subsequence, if necessary, we assume that $ f_{n_{0}}(x_n) = \| x_n \| = 1 $ for all $ n. $ Since
	\[ | f_{n_{0}} (x) - 1 | = | f_{n_{0}} (x) - f_{n_{0}} (x_n) | \leq \| f_{n_{0}} \| \| x_n - x \| \longrightarrow 0, \]
	we conclude that $ f_{n_{0}} (x) = \| x \| = 1, $ or, equivalently, $ f_{n_{0}} \in J(x). $ As $ \rho'_{-} (x, y) > 0, $ it is clear that $ f_{n_{0}} (y) = \rho'_{-} (x_{n_{0}}, y) > 0, $ as desired.\\
	Case II:  $ \{ f_n : n \in \mathbb{N} \} $ is infinite. Since $ B_{X^{*}} $ is weak*-compact, the set $ \{ f_n : n \in \mathbb{N} \} $ has an accumulation point, say, $ f_0. $ We claim that $ f_0 \in J(x). $ Let $ \epsilon > 0 $ be arbitrary. Since $ x_n \longrightarrow x, $ we can assume without any loss of generality that $ \| x - x_n \| < \epsilon $ for all $ n \in \mathbb{N}. $ Consider the weak*-open neighborhood $ \mathcal{O} := \{ f \in X^* : | \Psi(x) (f) - \Psi(x) (f_0) | < \epsilon \} \bigcap B_{X^{*}}, $ where $ \Psi:X \longrightarrow X^{**} $ is the canonical isometric embedding. Since $ f_0 $ is an accumulation point of the set $ \{ f_n : n \in \mathbb{N} \}, $ there exists $ m_0 \in \mathbb{N} $ such that $ f_{m_{0}} \in \mathcal{O}, $ and therefore, $ | f_{m_{0}} (x) - f_{0} (x) | < \epsilon. $ We also note that as before, $ | f_{m_{0}} (x) - 1 | = | f_{m_{0}} (x) - f_{m_{0}} (x_{m_{0}}) | \leq \| f_{m_{0}} \| \| x_{m_{0}} - x \| < \epsilon, $ which shows that $ f_{m_{0}} (x) > 1 - \epsilon. $ It is now immediate that $ f_0 (x) > f_{m_{0}} (x) - \epsilon > 1 - 2 \epsilon. $ The arbitrariness of $ \epsilon > 0 $ shows that $ f_0(x) \geq \| x \| = 1. $ Since $ f_0(x) \leq \| f_0 \| \| x \| = 1,$ it follows that $ f_0 \in J(x), $ thus justifying our claim.\\
	Let us now consider the weak*-open neighborhood $ \mathcal{O'} := \{ f \in X^* : | \Psi(x) (f) - \Psi(x) (f_0) | < \epsilon_0 \} \bigcap B_{X^{*}}, $ where $ 0 < \epsilon_0 < \frac{1}{2} \rho'_{-} (x, y). $ As before, it is easy to see that there exists $ n_0 \in \mathbb{N} $ such that $ | f_{n_{0}} (y) - f_0 (y) | < \epsilon_0, $ which implies that $ f_{n_{0}} (y) > f_0 (y) - \epsilon_0. $ Since $ \rho'_{-} (x, y) > 0 $ and $ f_0 \in J(x), $ we obtain that $ f_0 (y) \geq \rho'_{-} (x, y) > \epsilon_0.  $ Therefore, $ f_{n_{0}} (y) > f_0 (y) - \epsilon_0 > 0. $ As each $ x_n $ is smooth in $ X, $ it is clear that $ \rho'_{-} (x_{n_{0}}, y) = f_{n_{0}} (y) > 0, $ as desired. This establishes the theorem.
\end{proof}	

\begin{cor}\label{corollary:corollary4}(Proposition $ 2.7 $ of \cite{CCS})
	Suppose that $ X $ is a separable Banach space, and $ I $ is an index set with $ m $ elements. If there exists a ball-covering of $ S_{X} $ consisting of $ m $ balls, then there is a ball-covering $ \mathcal{B} = \{ B(x_i, r_i) : i \in I \} $ of $ S_{X} $ such that $ \{ x_i \}_{i \in I} $ are smooth points in $ X. $
\end{cor}

\begin{proof}
	The proof follows directly from Theorem \ref{theorem:smoothness}, by applying the well-known fact that smooth points are dense in a separable Banach space.
\end{proof}

We end this article with the following remark that further justifies the local approach to the ball-covering problem, that we have considered in this article:

\begin{remark}
	A detailed look at Proposition $ 2.7 $ of \cite{CCS} reveals that separability of the Banach space $ X $ is essential for the proof of the same. In contrast to that, Theorem \ref{theorem:smoothness} is essentially local in nature and the proof only requires the strictly weaker assumption of being able to approximate the concerned element by a sequence of smooth points. In particular, this shows that Theorem \ref{theorem:smoothness} is not only a theoretical generalization of Proposition $ 2.7 $ of \cite{CCS}, but it also allows us to tackle the ball-covering problem in non-separable Banach spaces, whenever the concerned element can be approximated by smooth points in the space. As a concrete example, let us consider the non-separable Banach space $ \ell_{\infty} $ and the non-smooth point $ x = (1, 1, ...) \in B_{\ell_{\infty}}. $ We observe that $ x $ can be approximated by smooth points in $ \ell_{\infty}. $ Indeed, let $ \{ u_n \} $ be the sequence of unit vectors in $ \ell_{\infty}, $ where for each $ n \in \mathbb{N}, $ $ u_n $ has $ 1 $ in the first coordinate and $ (1-\frac{1}{n}) $ in all other coordinates. Then it is easy to verify that each $ u_n $ is smooth in $ \ell_{\infty} $ and $ u_n \longrightarrow x $ as $ n \longrightarrow \infty. $ Therefore, Theorem \ref{theorem:smoothness} can be used effectively in this case to study ball-coverings, whereas Proposition $ 2.7 $ of \cite{CCS} is not applicable here. 
\end{remark}

	\bibliographystyle{amsplain}
	
\end{document}